\documentclass[11pt]{amsart}

\usepackage{amsmath}
\usepackage{amssymb}
\usepackage{amsfonts}
\usepackage{amsthm}
\usepackage{mathrsfs}
\usepackage{tikz-cd}
\usepackage[margin=1in]{geometry}
\usepackage{hyperref}  

\newtheorem{theorem}{Theorem}[section]
\newtheorem{proposition}[theorem]{Proposition}
\newtheorem{lemma}[theorem]{Lemma}
\newtheorem{corollary}[theorem]{Corollary}

\theoremstyle{definition}

\newtheorem{example}[theorem]{Example}

\theoremstyle{remark}
\newtheorem{remark}[theorem]{Remark}

\numberwithin{equation}{section}

\title{Chow groups of symmetric products of projective spaces}

\author{Youming Chen}
\address{School of Mathematical Science,
Chongqing University of Technology,
Chongqing 400054, P.R. China}
\email{youmingchen@cqut.edu.cn}

\author{Wenchuan Hu}
\address{School of Mathematics,
Sichuan University,
Chengdu 610064, P.R. China}
\email{huwenchuan@gmail.com}



\date{}

\begin{document}

\begin{abstract}
Let $\operatorname{SP}^{d}(\mathbb{P}^{n})$ denote the
$d$-th symmetric product of complex projective space.
We prove that its integral cycle class map is split injective
and has finite cokernel annihilated by $d!$.
We also compute the ranks of its Chow groups.
By constructing compatible integral retractions on Chow groups
and singular homology, we prove that addition of a fixed point
induces an isomorphism on the Chow group of $q$-dimensional cycles
whenever $d\geq q$.
For $n\geq 3$ and $d\geq 6$, we construct a nonalgebraic class
of order two in codimension two.
Its image generates a direct summand of the cycle class cokernel
and remains nonalgebraic after addition of fixed points.
Consequently, integral surjectivity fails both at certain finite
stages with $q\leq d$ and in the infinite symmetric product.
We further identify the codimension two cokernel with third
unramified cohomology of the smooth configuration space. 
Finally, we provide an example where multiplication by an even
integer annihilates a nonzero Chow class of order two.
\end{abstract}

\maketitle

\section{Introduction and main results}

Let $\mathbb{P}^{n}$ denote complex projective space of
dimension $n$.
For a nonnegative integer $d$, let
\[
X_{n,d}
=
\operatorname{SP}^{d}(\mathbb{P}^{n})
=
(\mathbb{P}^{n})^{d}/\mathfrak{S}_{d},
\]
where the symmetric group $\mathfrak{S}_{d}$ acts by permuting the factors.
We adopt  the convention $\operatorname{SP}^{0}(V)=\operatorname{Spec}(\mathbb{C})$.
Symmetric products are the Chow varieties parametrizing effective zero-cycles of fixed degree.
Their topology is closely related to the Dold--Thom theorem,
whereas their integral Chow groups are less well understood.
In particular, the comparison between algebraic cycles and integral singular homology involves torsion phenomena 
which are not visible with rational coefficients.

For a projective variety $X$, let $\operatorname{CH}_{i}(X)$
denote the Chow group of algebraic cycles of complex
dimension $i$, modulo rational equivalence.
We write
\[
\operatorname{cl}_{X,i}\colon
\operatorname{CH}_{i}(X)
\longrightarrow
H_{2i}(X,\mathbb{Z})
\]
for the integral cycle class map, and define
\[
D_i(X)
=
\operatorname{coker}(\operatorname{cl}_{X,i}).
\]
A homology class is called algebraic if it lies in the
image of the cycle class map.

For an abelian group $A$, we denote its torsion subgroup by
$A_{\mathrm{tors}}$, and write
\[
\operatorname{rank}A
=
\dim_{\mathbb{Q}}
\bigl(A\otimes_{\mathbb{Z}}\mathbb{Q}\bigr).
\]
The rational cycle class map of
$\operatorname{SP}^{d}(\mathbb{P}^{n})$ is an isomorphism.
This follows from the corresponding statement for
$(\mathbb{P}^{n})^{d}$ and the rational comparison for finite quotients; see \cite[Lem.~3.6]{ChenHu}.
The integral comparison was formulated in \cite[Conj.~1.2]{ChenHu}.
The discussion following that problem also proposed the
range $0\leq q\leq d$ as a candidate range for integral surjectivity.

Our first main result establishes injectivity and computes the ranks of the integral cycle class maps.

\begin{theorem}\label{theorem:main-comparison}
For all $n,d,i\geq 0$, the integral cycle class map
\[
\operatorname{cl}_{X_{n,d},i}\colon
\operatorname{CH}_{i}(X_{n,d})
\longrightarrow
H_{2i}(X_{n,d},\mathbb{Z})
\]
is split injective.
Its cokernel is finite and satisfies
$d!\,D_i(X_{n,d})=0.$
Consequently,
\[
H_{2i}(X_{n,d},\mathbb{Z})
\cong
\operatorname{CH}_{i}(X_{n,d})
\oplus D_i(X_{n,d}),
\]
and the cycle class map becomes an isomorphism after
tensoring with $\mathbb{Z}[1/d!]$.
The Chow groups are finitely generated, with ranks
\begin{equation}\label{equation:rank-formula}
r_{n,d}(i)
=
\#\left\{
(j_1,\ldots,j_d):
0\leq j_1\leq\cdots\leq j_d\leq n,\ 
\sum_{a=1}^{d}j_a=i
\right\}.
\end{equation}
\end{theorem}

The splitting in this theorem is a splitting of abelian groups.
We do not claim that it is canonical or compatible with any multiplicative structure.
For $d=0$, the empty tuple has sum zero, and $0!=1$.

The integral cycle class map is generally not surjective.

\begin{theorem}\label{theorem:main-obstruction}
Let $n\geq 3$, $d\geq 6$, and $N=nd$.
There exists a nonalgebraic class
\[
\alpha_{n,d}
\in H_{2N-4}(X_{n,d},\mathbb{Z})
\]
of exact order two whose image in $D_{N-2}(X_{n,d})$
generates a direct summand isomorphic to $\mathbb{Z}/2$.
In particular,
\begin{equation}\label{eq:decomp}
H_{2N-4}(X_{n,d},\mathbb{Z})
\cong
\operatorname{CH}_{N-2}(X_{n,d})
\oplus\mathbb{Z}/2\oplus D'_{n,d},
\end{equation}
where $D'_{n,d}$ is finite and annihilated by $d!$.
Both $\operatorname{CH}_{N-2}(X_{n,d})$ and
$H_{2N-4}(X_{n,d},\mathbb{Z})$ have rank two, 
yet they are not isomorphic even as abelian groups.
More precisely,
\[
\#H_{2N-4}(X_{n,d},\mathbb{Z})_{\mathrm{tors}}
\geq
2\,\#
\operatorname{CH}_{N-2}(X_{n,d})_{\mathrm{tors}}.
\]
\end{theorem}

The proof of Theorem~\ref{theorem:main-comparison} combines
Totaro's product and cyclic product lemmas
\cite[Lem. 6.2 and 8.1]{Totaro} with transfers for coarse
finite quotients.
These arguments also apply to projective varieties with
an affine paving.
The obstruction in
Theorem~\ref{theorem:main-obstruction} is detected by a
Steenrod square on the smooth configuration space.
The same construction works for each smooth connected
projective variety of dimension at least three.

\iffalse
Taking $(n,d)=(3,6)$ yields
\[
\operatorname{CH}_{16}
\bigl(\operatorname{SP}^{6}(\mathbb{P}^{3})\bigr)
\not\cong
H_{32}
\bigl(\operatorname{SP}^{6}(\mathbb{P}^{3}),\mathbb{Z}\bigr).
\]
Thus the unrestricted integral comparison in
\cite[Conj. 1.2]{ChenHu} is false.
\fi

The preceding example has cycle dimension $16>6$.
To study the proposed range $q\leq d$, we consider addition
of a fixed point.
Choose $o\in\mathbb{P}^{n}$ and let
\[
\iota_d\colon X_{n,d}\longrightarrow X_{n,d+1},
\qquad
c\longmapsto c+o.
\]

\begin{theorem}\label{theorem:main-stability}
For all $n,q\geq 0$ and $d\geq q$, addition of $o$ induces
an isomorphism
\[
(\iota_d)_*\colon
\operatorname{CH}_{q}(X_{n,d})
\xrightarrow{\ \cong\ }
\operatorname{CH}_{q}(X_{n,d+1}).
\]
The cycle class map on the infinite symmetric product is
split injective, with finite cokernel annihilated by $q!$
in cycle dimension $q$.
If $n\geq 3$, $m\geq 6$, and $q=nm-2$,
the class $\alpha_{n,m}$ in Theorem \ref{theorem:main-obstruction} remains nonalgebraic after addition
of any number of fixed points.
Its image generates a direct summand $\mathbb{Z}/2$ of $D_q(X_{n,d})$ for each $d\geq m$, and of the stable
cycle class cokernel.
 In particular,
\[
\operatorname{CH}_{16}\bigl(\operatorname{SP}^{16}(\mathbb{P}^3)\bigr)
 \not\cong H_{32}\bigl(\operatorname{SP}^{16}(\mathbb{P}^3),\mathbb{Z}\bigr),
\]
and the cycle class map on $\operatorname{SP}^\infty(\mathbb{P}^3)$ is not surjective in
homological degree $32$.
\end{theorem}

The core of the persistence argument is the construction of integral retractions 
on both Chow groups and singular homology which commute with the cycle class maps.
Topological split injectivity alone would not guarantee that a nonalgebraic class remains nonalgebraic.

We also compare the codimension-two cokernel with
unramified cohomology, following the argument of
Colliot-Thélène and Voisin
\cite[Thm.~3.7]{ColliotTheleneVoisin}.
The resulting classes on the configuration space ramify along the boundary of a smooth projective resolution.
Indeed, symmetric products of projective space are rational, 
and their smooth projective resolutions have vanishing integral Hodge defect in degree four.

Finally, multiplication of effective zero-cycles behaves
differently from addition of a fixed point.
We construct a nonzero class of order two in $\operatorname{CH}_{3}(\operatorname{SP}^{d}(\mathbb{P}^{2}))$ 
that is annihilated by each even multiplicity.
In particular, the doubling map
\[
\operatorname{SP}^{3}(\mathbb{P}^{2})
\longrightarrow
\operatorname{SP}^{6}(\mathbb{P}^{2})
\]
does not induce an injection on $\operatorname{CH}_{3}$.

\subsection*{Acknowledgements}

The authors thank Burt Totaro for early discussions on the
main topics of this paper and for suggesting arguments used
in the work.
This work is partially supported by the National Natural
Science Foundation of China
(Nos.~12126309 and 12126354), and by the Science and
Technology Research Program of Chongqing Municipal
Education Commission (No.~KJZD-K202601106).

\subsection*{Declaration concerning AI assistance}

ChatGPT was used as an auxiliary tool in manuscript
preparation.
The authors are responsible for the final text and
mathematical arguments.

\section{Integral transfers for finite quotients}\label{sec:transfers}

This section recalls basic facts about Chow groups and the cycle class maps.
In particular, we provide relative algebraic transfers for the Chow groups of the coarse quotient of a projective variety by a finite group,
along with compatible equivariant cell structures.
We also fix the conventions and notation for the arguments in what follows.
All varieties in this paper are reduced separated schemes
of finite type over $\mathbb{C}$.
They are not assumed to be irreducible.
Quotients by finite groups are coarse quotients.

\subsection{Cycle classes and localization}

For a variety $W$, let $\mathcal{Z}_i(W)$ denote the free abelian group generated by the integral closed subvarieties of $W$ of dimension $i$.
Then $\operatorname{CH}_i(W)$ is the quotient of $\mathcal{Z}_i(W)$ by rational equivalence.
For an integral variety $V$, let $\mathbb{C}(V)$ denote its
function field, and let
$\operatorname{div}_{V}(f)$ denote the principal divisor of
$f\in\mathbb{C}(V)^{\times}$.

For a possibly nonproper variety $W$, the cycle class map
takes values in Borel--Moore homology:
\[
\operatorname{cl}_{W,i}\colon
\operatorname{CH}_i(W)
\longrightarrow
H^{\mathrm{BM}}_{2i}(W,\mathbb{Z}).
\]
It commutes with proper pushforward and restriction to
open subsets; see \cite[Sec. 19.1]{Fulton}.
For proper $W$, Borel--Moore homology coincides with singular homology.
If $W$ is smooth of pure dimension $m$, we write
$\operatorname{CH}^{r}(W)=\operatorname{CH}_{m-r}(W)$
and apply Poincar\'e duality to define
\[
\operatorname{cl}^{r}_{W}\colon
\operatorname{CH}^{r}(W)
\longrightarrow
H^{2r}(W,\mathbb{Z}).
\]
Here Poincar\'e duality will only be applied to smooth varieties.

Let $B\subset W$ be closed, where $W$ has pure dimension
$m$, and suppose that $\dim_{\mathbb{C}}B\leq m-c$.
Set $U=W\setminus B$.
The localization sequence for Chow groups yields
\begin{equation}\label{eq:loc-chow}
\operatorname{CH}_{m-r}(W)
\xrightarrow{\ \cong\ }
\operatorname{CH}_{m-r}(U)
\qquad (r<c).
\end{equation}
Meanwhile, the localization sequence for Borel--Moore homology gives
\begin{equation}\label{eq:loc-bm}
H^{\mathrm{BM}}_{k}(W,A)
\xrightarrow{\ \cong\ }
H^{\mathrm{BM}}_{k}(U,A)
\qquad
\bigl(k>2\dim_{\mathbb{C}}B+1\bigr)
\end{equation}
for each coefficient group $A$.
Indeed, the adjacent homology groups of $B$ vanish by
the dimension bound.
If $W$ is smooth, this implies
\[
H^{s}(W,A)\xrightarrow{\ \cong\ }H^{s}(U,A)
\qquad (s\leq 2c-2).
\]

\subsection{Transfers on Chow groups}

Let a finite group $G$ act algebraically on a projective
variety $M$.
For $K\leq G$, write
\[
\pi_K\colon M\longrightarrow M/K.
\]
If $W\subset M$ is integral, set
\[
S_K(W)=\{g\in K:gW=W\},
\qquad
I_K(W)=\{g\in K:g|_W=\operatorname{id}_W\},
\]
and
\[
e_K(W)=|I_K(W)|.
\]
The action need not be faithful.

If $V\subset M/K$ is integral, the reduced irreducible
components of $\pi_K^{-1}(V)$ form one $K$-orbit and dominate $V$.
For any such component $W$,
\begin{equation}\label{equation:generic-degree}
[\mathbb{C}(W):\mathbb{C}(V)]
=
\frac{|S_K(W)|}{e_K(W)}.
\end{equation}
This follows by considering a general fiber and its
stabilizer.
In particular, we have the following result on the weighted inverse images.

\begin{lemma}%\label{lemma:weighted-pullback}
The homomorphism on cycle groups defined by
\[
\pi_K^{\#}[V]
=
\sum_{W\subset\pi_K^{-1}(V)}
e_K(W)[W]
\]
preserves rational equivalence.
Here and below, such sums runs over reduced irreducible components, 
and the symbol $\#$ denotes this weighted inverse image, not flat pullback.
Moreover,
\[
(\pi_K)_*\pi_K^{\#}=|K|\operatorname{id},
\qquad
\pi_K^{\#}(\pi_K)_*
=
\sum_{g\in K}g_*
\]
on cycle groups and hence on Chow groups.
\end{lemma}

\begin{proof}
The first identity follows from
\eqref{equation:generic-degree} and the number
$\frac{|K|}{|S_K(W)|}$ of components above $V$.
For the second identity, each distinct translate of $W$
has coefficient $|S_K(W)|$ on both sides.
Observe the pushforward  on the invariant cycle group
\[
(\pi_K)_*\colon Z_i(M)^K\longrightarrow Z_i(M/K)
\]
is injective.
Indeed, an orbit sum maps to a positive multiple of a
basis element, and distinct orbit sums map to distinct
basis elements.
Let $V\subset M/K$ be integral of dimension $i+1$, and let
$f\in\mathbb{C}(V)^{\times}$.
We claim that
\begin{equation}\label{equation:principal-transfer}
\pi_K^{\#}\operatorname{div}_{V}(f)
=
\sum_{W\subset\pi_K^{-1}(V)}
e_K(W)\operatorname{div}_{W}(f\circ\pi_K).
\end{equation}
Both sides are invariant cycles.
Their pushforwards are equal to
$|K|\operatorname{div}_{V}(f)$ by the norm formula for
principal divisors
\cite[Sec. 1.4]{Fulton}.
Injectivity on invariant cycle groups proves
\eqref{equation:principal-transfer}.
Its right-hand side is rationally equivalent to zero.
The argument also applies to nonnormal subvarieties,
using the usual definition of principal divisors through
normalization.
\end{proof}

\begin{proposition}%\label{proposition:relative-transfer}
Let $H\leq G$, and let $q\colon M/H\longrightarrow M/G$.
There exists a homomorphism
\[
\operatorname{tr}^{G}_{H}\colon
\operatorname{CH}_i(M/G)
\longrightarrow
\operatorname{CH}_i(M/H)
\]
such that
\[
q_*\operatorname{tr}^{G}_{H}
=
[G:H]\operatorname{id}.
\]
On cycle groups it is given by
\begin{equation}\label{equation:relative-transfer}
\operatorname{tr}^{G}_{H}[V]
=
\sum_{U\subset q^{-1}(V)}
\frac{e_G(W_U)}{e_H(W_U)}[U],
\end{equation}
where $W_U$ is any irreducible component of
$\pi_H^{-1}(U)$.
\end{proposition}

\begin{proof}
Since $I_H(W_U)=H\cap I_G(W_U)$, all coefficients are integers.
They are independent of the choice of $W_U$ within its
$H$-orbit.
The formula is characterized on cycle groups by
\begin{equation}\label{equation:transfer-characterization}
\pi_H^{\#}\operatorname{tr}^{G}_{H}
=
\pi_G^{\#}.
\end{equation}
The map $\pi_H^{\#}$ is injective on cycle groups.
For an integral subvariety $V$ and a rational function
$f$ on $V$, write
$a_U=\frac{e_G(W_U)}{e_H(W_U)}$.
Applying $\pi_H^{\#}$ and using
\eqref{equation:principal-transfer} gives
\[
\operatorname{tr}^{G}_{H}\operatorname{div}_{V}(f)
=
\sum_{U\subset q^{-1}(V)}
a_U\operatorname{div}_{U}(f\circ q).
\]
Thus the map preserves rational equivalence.

Applying $(\pi_G)_*$ to
\eqref{equation:transfer-characterization} gives
\[
|H|q_*\operatorname{tr}^{G}_{H}
=
|G|\operatorname{id}
\]
on the free cycle group.
Cancellation there proves the stated identity.
\end{proof}

For a free action, this transfer agrees with flat pullback
for the finite \'etale map $q$.
In general, $q$ need not be flat.
The stabilizer multiplicities are necessary in general.
For example, when $M$ is a point, the transfer is
multiplication by $[G:H]$, although $q$ is the identity
map of the underlying variety.

\subsection{Transfers on homology}

We apply finite semialgebraic equivariant cell structures.
The relevant construction and subdivision results are
given in \cite[Prop. 2.2 and Lem. 2.6]{ParkSuh}.

\begin{lemma}\label{lemma:equivariant-cells}
Let a finite group $G$ act on a compact semialgebraic
subset $M$ of an orthogonal real representation.
Given finitely many closed semialgebraic subsets of $M$,
there is a finite semialgebraic $G$-CW structure for which
each of them is a subcomplex of the underlying CW complex.
Every subgroup quotient has the induced finite CW
structure.
Any two such structures admit a common equivariant
subdivision, and subdivision induces the identity on the
homology of each subgroup quotient.
\end{lemma}

\begin{proof}
Adjoin all translates of the prescribed subsets and write
the resulting family as $A_1,\ldots,A_s$.
For each $r$, let
\[
B_r
=
\{x\in M:\#\{j:x\in A_j\}=r\}.
\]
These are invariant semialgebraic subsets.
The intersection $A_j\cap B_r$ is both closed and open in
$B_r$.
Thus membership in each $A_j$ is constant on each
connected subset of $B_r$.

In the construction of
\cite[Prop. 2.2]{ParkSuh}, choose the quotient
triangulation compatibly with the images of the $B_r$
and with the orbit types.
Every underlying open cell is connected and lies in one
$B_r$.
It therefore lies in $A_j$ or is disjoint from $A_j$.
Closedness of $A_j$ makes it a subcomplex.

For a finite group $G$, an equivariant open cell has the form
$G/S\times\operatorname{int}D^k$, with $S$ acting
trivially on the disk.
This yields the quotient cell structures.
The common subdivision statement follows from
\cite[Lem. 2.6]{ParkSuh}.
\end{proof}

The lemma may apply to projective $G$-varieties.
Indeed, a compact semialgebraic embedding
$e\colon M\longrightarrow\mathbb{R}^{a}$ yields an
equivariant embedding
\[
x\longmapsto
\bigl(e(g^{-1}x)\bigr)_{g\in G}
\]
into a real permutation representation.

Let $C_{\bullet}=C^{\mathrm{cell}}_{\bullet}(M,\mathbb{Z})$  for a structure as in Lemma~\ref{lemma:equivariant-cells}.
For a $K$-module $A$, let
\[
A_K=A/\langle ga-a:g\in K,\ a\in A\rangle
\]
be its group of coinvariants.
The quotient cell structures give
\begin{equation}\label{equation:cell-coinvariants}
C^{\mathrm{cell}}_{\bullet}(M/K,\mathbb{Z})
\cong (C_{\bullet})_K.
\end{equation}
Orientations are chosen on the disk factor of each equivariant cell, so the identifications respect
the boundary maps.
\begin{proposition}%\label{proposition:homological-transfer}
For each $H\leq G$, there are homomorphisms
\[
\operatorname{Tr}^{G}_{H}\colon
H_k(M/G,\mathbb{Z})
\longrightarrow H_k(M/H,\mathbb{Z})
\]
which are independent of the chosen cell structure and
satisfy
\[
q_*\operatorname{Tr}^{G}_{H}
=
[G:H]\operatorname{id}.
\]
They are compatible with cycle classes:
\[
\operatorname{cl}_{M/H,i}
\operatorname{tr}^{G}_{H}
=
\operatorname{Tr}^{G}_{H}
\operatorname{cl}_{M/G,i}.
\]
\end{proposition}

\begin{proof}
On the coinvariant chain complexes,  define the coset sum
\begin{equation}\label{equation:coset-transfer}
 (C_\bullet)_G\longrightarrow(C_\bullet)_H,\qquad
 [c]_G\longmapsto\sum_{Hg\in H\backslash G}[gc]_H.
\end{equation}
Changing representatives does not affect the sum.
Replacing $c$ by $ac$, with $a\in G$, permutes the cosets.
The map is therefore well defined and commutes with the
boundary.
Its composite with the quotient map is multiplication by
$[G:H]$.
Equivariant subdivisions commute with
\eqref{equation:coset-transfer}, so
Lemma~\ref{lemma:equivariant-cells} yields independence of
the cell structure.

To verify compatibility with cycle classes, let
$V\subset M/G$ be integral of dimension $i$.
Choose a compatible cell structure after removing from
the components above $V$ their singular loci,
intersections, loci where stabilizers increase, and the
exceptional loci of the generically \'etale maps.
These subsets have real dimension at most $2i-2$.

If a lift of a top-dimensional cell has stabilizer $S$,
its images in $M/H$ are indexed by double cosets $HgS$.
The number of terms in
\eqref{equation:coset-transfer} contributing to the
corresponding cell is
\[
[S:S\cap g^{-1}Hg]
=
\frac{|gSg^{-1}|}{|H\cap gSg^{-1}|}.
\]
At a general point of the component above $U$, this is
$\frac{e_G(W_U)}{e_H(W_U)}$.
The maps preserve complex orientations.
Thus the transfer of the fundamental class of $V$
is the cycle class of
\eqref{equation:relative-transfer}.
\end{proof}

The chain identification
\eqref{equation:cell-coinvariants} does not imply an
integral identification of quotient homology with
coinvariants of homology.
Taking coinvariants is exact over $\mathbb{Q}$, but not
in general over $\mathbb{Z}$.

\subsection{Proper equivariant maps}
The following compatibility will be essential when a configuration acquires
an additional fixed point. In particular, it applies to maps whose images
meet collision strata.

\begin{proposition}\label{proposition:transfer-naturality}
Let $f\colon M\longrightarrow M'$ be a proper
$G$-equivariant morphism of projective $G$-varieties.
For $K\leq G$, let
$f_K\colon M/K\longrightarrow M'/K$ be the induced map.
Then, for each $H\leq G$,
\[
(f_H)_*\operatorname{tr}^{G}_{H}
=
\operatorname{tr}^{G}_{H}(f_G)_*,
\qquad
(f_H)_*\operatorname{Tr}^{G}_{H}
=
\operatorname{Tr}^{G}_{H}(f_G)_*.
\]
\end{proposition}

\begin{proof}
For a cycle $z$ on $M/K$, the cycles $f_*\pi_{M,K}^{\#}(z)$ and $\pi_{M',K}^{\#}(f_K)_*(z)$ are $K$-invariant.
Their pushforwards to $M'/K$ are both
$|K|(f_K)_*(z)$.
Injectivity of pushforward on invariant cycle groups
therefore yields
\[
f_*\pi_{M,K}^{\#}
=
\pi_{M',K}^{\#}(f_K)_*.
\]
Apply this with $K=H$ and $K=G$, and apply
\eqref{equation:transfer-characterization}.
Injectivity of $\pi_{M',H}^{\#}$ proves the first
identity.
The same argument includes cycles contracted by $f$.

For homology, choose finite $G$-CW structures and replace
the underlying continuous map of $f$ by an equivariantly
homotopic cellular map, using equivariant cellular
approximation.
The induced maps on each subgroup quotient are
homotopic to the original quotient maps.
The cellular chain map is $G$-linear, and each
$G$-linear chain map commutes with the coset sum
\eqref{equation:coset-transfer}.
This proves the second identity.
\end{proof}

\section{Split injectivity and ranks}\label{sec:splitting}

In this section we prove the split injectivity of Theorem~\ref{theorem:main-comparison}.
First we give a Sylow descent result for finite groups.

\subsection{Sylow descent}

\begin{proposition}\label{proposition:sylow-descent}
Let a nontrivial finite group $G$ act on a projective
variety $M$.
For each prime $\ell$ dividing $|G|$, choose a Sylow
subgroup $P_{\ell}$.
If the cycle class map of $M/P_{\ell}$ is split
injective in dimension $i$ for each such $\ell$, then
the cycle class map of $M/G$ is split injective in that
dimension.
The same assertion holds for injectivity without a
splitting.
\end{proposition}

\begin{proof}
Set
\[
X=M/G,\qquad
Y_{\ell}=M/P_{\ell},\qquad
m_{\ell}=[G:P_{\ell}],
\]
and let $q_{\ell}\colon Y_{\ell}\longrightarrow X$.
The integers $m_{\ell}$ have greatest common divisor one.
Indeed, a common prime divisor would divide $|G|$,
whereas $\ell$ does not divide $m_{\ell}$.
Choose integers $b_{\ell}$ such that
$
\sum_{\ell\mid |G|}b_{\ell}m_{\ell}=1.
$
Let $s_{\ell,i}$ be a left inverse of the cycle class map
of $Y_{\ell}$.
Then
\[
S_i
=
\sum_{\ell\mid |G|}
b_{\ell}(q_{\ell})_*s_{\ell,i}
\operatorname{Tr}^{G}_{P_{\ell}}
\]
is a left inverse of $\operatorname{cl}_{X,i}$, by
compatibility of transfers with cycle classes:
\[
 S_i  \operatorname{cl}_{X,i}
 =\sum_\ell b_\ell q_{\ell\ast}s_{\ell,i}
       \operatorname{cl}_{Y_\ell,i}  \operatorname{tr}_{P_\ell}^G
 =\sum_\ell b_\ell m_\ell  \mathrm{id}=  \mathrm{id}.
\]
If only injectivity is assumed, a class in the kernel
transfers to zero on each $Y_{\ell}$.
It is therefore annihilated by each $m_{\ell}$ and hence vanishes.
\end{proof}

Let $\mathcal{A}$ denote the class of projective varieties 
which admit a finite partition into locally closed subvarieties 
isomorphic to Zariski open subsets of affine spaces and whose cycle class maps are split injective in each dimension.
For a prime $\ell$, set
\[
\operatorname{Cyc}_{\ell}(V)
=
V^{\ell}/(\mathbb{Z}/\ell),
\]
where the action is cyclic permutation.

\begin{lemma}%\label{lemma:totaro-closure}
The class $\mathcal{A}$ is closed under finite products
and under the construction $V\longmapsto\operatorname{Cyc}_{\ell}(V)$.
\end{lemma}

\begin{proof}
The product assertion follows from
\cite[Lem.~6.2]{Totaro}.
For cyclic products,
\cite[Lem.~8.1(2)]{Totaro} yields the required partition
together with a surjection from its auxiliary integral graded group onto
$\operatorname{CH}_{*}(\operatorname{Cyc}_{\ell}(V))$.
By part~(3) of that lemma, the composite with the cycle
class map is split injective.
The surjection is therefore an isomorphism, and the
cycle class map is split injective.
\end{proof}

\begin{theorem}\label{theorem:affine-paving-comparison}
If $V\in\mathcal{A}$, each cycle class map of
$\operatorname{SP}^{d}(V)$ is split injective.
In particular, this holds for each projective variety $V$ admitting an affine paving.
\end{theorem}

\begin{proof}
The cases $d=0,1$ are immediate.
For a prime $\ell$, define
\[
W_{\ell,0}=\{1\},
\qquad
W_{\ell,r+1}
=
W_{\ell,r}^{\ell}\rtimes(\mathbb{Z}/\ell),
\]
where the cyclic group permutes the $\ell$ blocks.
Also define
\[
V_{\ell,0}=V,
\qquad
V_{\ell,r+1}
=
\operatorname{Cyc}_{\ell}(V_{\ell,r}).
\]
Then quotienting first by the normal subgroup yields
\[
V^{\ell^{r}}/W_{\ell,r}
\cong V_{\ell,r},
\]
and each $V_{\ell,r}$ lies in $\mathcal{A}$.
Write
$d=\sum_{r\geq 0}a_r\ell^r, 0\leq a_r<\ell.$
The subgroup
\[
P_{\ell}
=
\prod_{r\geq 0}W_{\ell,r}^{a_r}
\subset\mathfrak{S}_{d}
\]
acts on disjoint blocks and forms a Sylow subgroup.
Indeed,
\[
|W_{\ell,r}|
=
\ell^{(\ell^r-1)/(\ell-1)}, \quad
v_{\ell}(|P_{\ell}|)
=
\frac{d-\sum_r a_r}{\ell-1}
=
\sum_{j\geq 1}\left\lfloor\frac{d}{\ell^j}\right\rfloor
=
v_{\ell}(d!).
\]
Moreover,
\[
V^d/P_{\ell}
\cong
\prod_{r\geq 0}V_{\ell,r}^{a_r}
\in\mathcal{A}.
\]
Proposition~\ref{proposition:sylow-descent} establishes the first assertion.

If $V$ admits an affine paving, localization along the paving
shows that its odd homology vanishes and its even homology
is free.
Comparison with Chow localization shows that its cycle
class maps are isomorphisms.
The affine strata give the required partition, so $V\in\mathcal{A}$.
\end{proof}

\subsection{Factorial bounds and rank formulas}

\begin{proposition}\label{proposition:bounds-ranks}
For $X=X_{n,d}$, the group $D_i(X)$ is finite and
annihilated by $d!$.
All torsion in $\operatorname{CH}_i(X)$ and
$H_k(X,\mathbb{Z})$ is annihilated by $d!$, and each odd homology group is finite.
The ranks are given by
\eqref{equation:rank-formula}, or equivalently by
\[
\sum_{d,i\geq 0}r_{n,d}(i)t^d u^i
=
\prod_{j=0}^{n}(1-tu^j)^{-1}.
\]
\end{proposition}

\begin{proof}
Let $M=(\mathbb{P}^{n})^d$, and let
$\pi\colon M\longrightarrow X$.
The product affine paving identifies
$\operatorname{CH}_i(M)$ with the free group
$H_{2i}(M,\mathbb{Z})$, and the odd homology of $M$
vanishes.
For $z\in H_{2i}(X,\mathbb{Z})$, choose
$\zeta\in\operatorname{CH}_i(M)$ whose cycle class equals $\operatorname{Tr}^{\mathfrak{S}_d}_{\{1\}}(z)$.
Then
\[
\operatorname{cl}_{X,i}(\pi_*(\zeta))=\pi_\ast\operatorname{Tr}_{\{1\}}^G(z)=d!z.
\]
Thus $d!D_i(X)=0$.
Since $X$ is a compact semialgebraic space, its homology
is finitely generated.
Therefore, the cokernel is finite.

Torsion classes transfer to zero inside the corresponding torsion-free groups of $M$.
Pushing forward shows that they are annihilated by $d!$.
Hence, odd homology also transfers to zero and is finite.

By Theorem~\ref{theorem:affine-paving-comparison},
$\operatorname{CH}_i(X)$ is a direct summand of
$H_{2i}(X,\mathbb{Z})$.
Its rank is the rational even Betti number.
Taking coinvariants is exact over $\mathbb{Q}$, so
\[
H_*(X,\mathbb{Q})
\cong
H_*(M,\mathbb{Q})_{\mathfrak{S}_d}.
\]
The Künneth basis is indexed by tuples $(j_1,\ldots,j_d)$ with $0\leq j_a\leq n$.
All degrees are even, so permutations introduce no
signs.
The orbit representatives are precisely the nondecreasing tuples appearing in \eqref{equation:rank-formula}.
Their multiplicities yield the stated generating series.
\end{proof}

\begin{proof}[Proof of Theorem~\ref{theorem:main-comparison}]
Apply Theorem~\ref{theorem:affine-paving-comparison} to $\mathbb{P}^{n}$ and invoke Proposition~\ref{proposition:bounds-ranks}.
For $i>nd$, both groups vanish.
\end{proof}

The same proof gives the following extension.

\begin{corollary}\label{corollary:paving-bounds}
Let $V$ be projective and admit an affine paving.
For $\operatorname{SP}^{d}(V)$, each cycle class
cokernel is finite and annihilated by $d!$.
All Chow torsion and homology torsion are annihilated by
$d!$, and the odd homology groups are finite.

If $b_j=\operatorname{rank}H_{2j}(V,\mathbb{Z})$, then
\[
\sum_{d,i\geq 0}
\operatorname{rank}
\operatorname{CH}_i(\operatorname{SP}^{d}(V))
t^d u^i
=
\prod_{j\geq 0}(1-tu^j)^{-b_j}.
\]
\end{corollary}

These statements do not imply that the Chow groups are
torsion free.

\begin{example}
The Sylow subgroups of $\mathfrak{S}_6$ for primes $2,3,5$ have indices $45,80,144$, respectively.
Since
\[
5\cdot45-80-144=1,
\]
a left inverse in
Proposition~\ref{proposition:sylow-descent} is
\[
S_i
=
5(q_2)_*s_{2,i}\operatorname{Tr}^{\mathfrak{S}_6}_{P_2}
-
(q_3)_*s_{3,i}\operatorname{Tr}^{\mathfrak{S}_6}_{P_3}
-
(q_5)_*s_{5,i}\operatorname{Tr}^{\mathfrak{S}_6}_{P_5}.
\]
All coefficients and maps are integral.
\end{example}

\section{Nonalgebraic torsion in codimension two}\label{sec:obstruction}

This section constructs a nonalgebraic homology class in the singular homology of the symmetric product of projective space,
which disproves surjectivity of the cycle class map.

\subsection{A Steenrod square obstruction}

For a real vector bundle $E$, let
$w_j(E)\in H^j(-,\mathbb{F}_2)$ denote its
Stiefel--Whitney classes.
For a complex vector bundle, let $c_j$ denote its Chern
classes.
Let
\[
\rho_2\colon H^*(B,\mathbb{Z})
\longrightarrow H^*(B,\mathbb{F}_2)
\]
be reduction modulo two, and let
\[
\beta_2\colon H^a(B,\mathbb{F}_2)
\longrightarrow H^{a+1}(B,\mathbb{Z})
\]
be the connecting homomorphism associated to the short exact sequence
$0\longrightarrow\mathbb{Z}\xrightarrow{\ 2\ }\mathbb{Z}
\longrightarrow\mathbb{F}_2\longrightarrow 0.$
We use the Steenrod squares
$\operatorname{Sq}^s:H^a(B,\mathbb{F}_2)\to H^{a+s}(B,\mathbb{F}_2)$ together with the identities
\begin{equation}\label{eq:squares}
 \rho_2\beta_2=\operatorname{Sq}^1,\qquad
\operatorname{Sq}^s(ab)=\sum_{r=0}^s\operatorname{Sq}^r(a)\operatorname{Sq}^{s-r}(b),\qquad
\operatorname{Sq}^s(x^a)=\binom as x^{a+s}\quad(|x|=1).
\end{equation}
These are the Bockstein identity, Cartan formula, and degree one
polynomial formula \cite[Sec. 3.E and 4.L]{HatcherTopology}.
In particular, $\operatorname{Sq}^3$ is the Steenrod square of index three, raising degree
by three; it is not an iteration of $\operatorname{Sq}^1$.

\begin{lemma}\label{lemma:steenrod-vanishing}
Let $M$ be a smooth complex quasiprojective variety.
For each $\zeta\in\operatorname{CH}^{2}(M)$,
\[
\operatorname{Sq}^{3}\rho_2
\operatorname{cl}^{2}_{M}(\zeta)=0.
\]
\end{lemma}

\begin{proof}
It suffices to consider an integral codimension-two
subvariety $Z\subset M$.
A finite locally free resolution of $\mathcal{O}_Z$
defines a virtual algebraic bundle
$v_Z=[\mathcal{O}_Z]\in K_0(M)$.
The formula for the first nonzero Chern class yields
\[
c_1(v_Z)=0,
\qquad
c_2(v_Z)=-[Z].
\]
See \cite[Ex. 15.3.6]{Fulton} and \cite[proof of Cor. 3.2]{Totaro}.
This formula does not require $Z$ to be smooth or a local complete intersection.
Hence, for the underlying virtual topological bundle,
\[
\operatorname{cl}^{2}_{M}([Z])
=
-c_2(v_Z^{\mathrm{top}}).
\]
Moreover, stable virtual complex bundles are classified by $\mathbb{Z}\times BU$.
Its mod-two cohomology is concentrated in even degrees;
see \cite[Thm. 3.9]{HatcherBundles}.
In particular, its degree-seven cohomology group vanishes.
Therefore, naturality gives the claimed vanishing.
\end{proof}

The codimension-two condition is essential here:
the coefficient $(r-1)!$ appearing in the first nonzero Chern class
formula equals one for $r=2$.
This lemma does not assert that all positive Steenrod operations vanish on arbitrary algebraic cycle classes.

\subsection{The configuration space}

Assume $n\geq 3$ and $d\geq 6$.
Set $N=nd$ and
\[
F_d(\mathbb{P}^{n})
=
(\mathbb{P}^{n})^d
\setminus
\bigcup_{r<s}\{p_r=p_s\},
\qquad
U_{n,d}
=
F_d(\mathbb{P}^{n})/\mathfrak{S}_d.
\]
Here $p_r$ denotes the $r$-th coordinate point.
Then the quotient map
$F_d(\mathbb{P}^{n})\longrightarrow U_{n,d}$ is finite \'etale of degree $d!$,
and $U_{n,d}$ is a smooth open subvariety of $X_{n,d}$ of dimension $N=nd$.
The collision locus $\Delta_X=X_{n,d}\setminus U_{n,d}$ has dimension $N-n$.
Since $n\geq3$, Localization formulas together with \eqref{eq:loc-chow} and \eqref{eq:loc-bm}  yield a commutative diagram
\begin{equation}\label{equation:configuration-comparison}
\begin{tikzcd}
\operatorname{CH}_{N-2}(X_{n,d})
  \arrow[r,"\operatorname{cl}_{X_{n,d},N-2}"]
  \arrow[d,"\cong"']
&
H_{2N-4}(X_{n,d},\mathbb{Z})
  \arrow[d,"\Psi","\cong"']
\\
\operatorname{CH}^{2}(U_{n,d})
  \arrow[r,"\operatorname{cl}^{2}_{U_{n,d}}"']
&
H^4(U_{n,d},\mathbb{Z}).
\end{tikzcd}
\end{equation}
Here $\Psi$ is restriction in Borel--Moore homology composed with inverse Poincar\'e duality.
Explicitly,
$2N-4>2(N-n)+1$, so both adjacent terms in Borel--Moore localization vanish. No duality on $X_{n,d}$ is invoked.
Moreover, the real permutation representation gives a bundle
\[
E
=
F_d(\mathbb{P}^{n})
\times_{\mathfrak{S}_d}\mathbb{R}^{d}
\longrightarrow U_{n,d}.
\]
Its fiber admits a basis indexed by the points of the unordered configuration. 
Set
\[
a=\beta_2(w_3(E))\in H^4(U_{n,d},\mathbb{Z}).
\]
The coefficient sequence implies $2a=0$.

\begin{proposition}\label{proposition:explicit-obstruction}
The class $a$ has exact order two and is not algebraic.
Consequently,
\[
\alpha_{n,d}=\Psi^{-1}(a)
\]
is a nonalgebraic class of exact order two  in $H_{2N-4}(X_{n,d},\mathbb{Z})$.
\end{proposition}

\begin{proof}
Choose three disjoint small Euclidean balls inside an affine chart $\mathbb{C}^{n}\subset\mathbb{P}^{n}$, 
with centers $c_1,c_2,c_3$.
Choose $\epsilon>0$ such that the radius $\epsilon$ spheres lie in the respective balls.
Choose $d-6$ distinct fixed points $z_7,\ldots,z_d$ outside their closures.
Let $T=(\mathbb{RP}^{2n-1})^3.$
For unit vectors $v_r\in\mathbb{C}^{n}$, define
\[
j([v_1],[v_2],[v_3])
=
\sum_{r=1}^{3}
\bigl([c_r+\epsilon v_r]+[c_r-\epsilon v_r]\bigr)
+
\sum_{s=7}^{d}[z_s].
\]
On the left-hand side, brackets denote a real line; on the right-hand side, they denote
a degree-one effective zero-cycle.
Changing the sign of $v_r$ swaps the two points in the corresponding pair.
Thus $j:T\to U_{n,d}$ is well defined and continuous.

Let $L_r$ be the pullback of the tautological real line bundle from the $r$-th factor.
The two-point permutation bundle of a moving pair decomposes as $\varepsilon^1\oplus L_r$, 
using the sum and difference of its basis vectors.
Therefore
\[
j^*E
\cong
\varepsilon^{d-3}\oplus L_1\oplus L_2\oplus L_3,
\]
where $\varepsilon^{d-3}$ denotes the trivial real bundle of rank $d-3$.
Write
\[
x=w_1(L_1),\qquad
y=w_1(L_2),\qquad
z=w_1(L_3).
\]
Then
\[
H^*(T,\mathbb{F}_2)
=
\mathbb{F}_2[x,y,z]/(x^{2n},y^{2n},z^{2n}),
\qquad |x|=|y|=|z|=1.
\]
By the Whitney sum formula, $j^\ast w_3(E)=xyz$, whence
\[
\rho_2j^*(a)
=\operatorname{Sq}^1(xyz)=
x^2yz+xy^2z+xyz^2.
\]
Applying the identities in \eqref{eq:squares}, we obtain 
\begin{align*}
 \operatorname{Sq}^3(x^2yz)&=x^4y^2z+x^4yz^2,\\
\operatorname{Sq}^3(xy^2z)&=x^2y^4z+xy^4z^2,\\
\operatorname{Sq}^3(xyz^2)&=x^2yz^4+xy^2z^4.
\end{align*}
For example, $\operatorname{Sq}^1(x^2)=0$, $\operatorname{Sq}^2(x^2)=x^4$, and
$\operatorname{Sq}^s(y)=\operatorname{Sq}^s(z)=0$ for $s>1$, 
so the first identity has exactly the two displayed terms. 
Thus
\begin{equation*}%\label{eq:detector}
 j^\ast\operatorname{Sq}^3\rho_2(a)= \operatorname{Sq}^3\rho_2(j^\ast(a))
 =x^4y^2z+x^4yz^2+x^2y^4z+xy^4z^2+x^2yz^4+xy^2z^4.
\end{equation*}
These six monomials are distinct, and each exponent is at most $4<2n$.
None is truncated or cancelled, so the expression is nonzero.
Lemma~\ref{lemma:steenrod-vanishing} implies its nonalgebraicity.
Since $a\neq0$ and $2a=0$, the order is exactly two.
The commutative diagram~\eqref{equation:configuration-comparison} then gives the claim for $\alpha_{n,d}$.
\end{proof}

\begin{proof}[Proof of Theorem~\ref{theorem:main-obstruction}]
Let
$\lambda\colon H^7(T,\mathbb{F}_2)\longrightarrow\mathbb{F}_2$
extract the coefficient of $x^4y^2z$ in the monomial basis.
Then  the homomorphism
\[
\chi=\lambda j^*\operatorname{Sq}^{3}\rho_2\Psi\colon H_{2N-4}(X_{n,d},\mathbb{Z}) \longrightarrow\mathbb{Z}/2
\]
vanishes on algebraic classes by the diagram~\eqref{equation:configuration-comparison} 
and Lemma~\ref{lemma:steenrod-vanishing}, while $\chi(\alpha_{n,d})=1$.
It descends to $\overline{\chi}\colon D_{N-2}(X_{n,d})\longrightarrow\mathbb{Z}/2.$

Since $\bar\alpha_{n,d}$ has order two and maps to $1$, the map
$1\mapsto\bar\alpha_{n,d}$ defines a section.
Set $D'_{n,d}=\ker\bar\chi$; this gives
\[
 D_{N-2}(X)\cong \mathbb{Z}/2\oplus D'_{n,d}.
\]
Theorem~\ref{theorem:main-comparison} provides finiteness, the factorial annihilator,
and \eqref{eq:decomp}.
The coefficient functional is essential: an order-two element in an
abelian group need not generate a direct summand.

The codimension-two tuples in \eqref{equation:rank-formula} have exactly two permutation
orbits, represented in codimension notation by
$(2,0,\ldots,0)$ and $(1,1,0,\ldots,0)$.
Both orbits occur, so the two groups have rank two.
Taking torsion in \eqref{eq:decomp} gives
\[
 \#H_{2N-4}(X,\mathbb{Z})_{\operatorname{tors}}
 =2\,\#\operatorname{CH}_{N-2}(X)_{\operatorname{tors}}\,\#D'_{n,d}.
\]
This establishes the strict difference in torsion cardinalities and hence
abstract nonisomorphism.
\end{proof}

\begin{remark}
The summand $\mathbb{Z}/2$ in the homology decomposition
may be generated by $\alpha_{n,d}$ itself.
If $R$ is a left inverse of the cycle class map, define
$s\colon\mathbb{Z}/2\longrightarrow
\operatorname{CH}_{N-2}(X_{n,d})$ by
$s(1)=R(\alpha_{n,d})$.
Then $R-s\chi$ is a left inverse that annihilates $\alpha_{n,d}$.
\end{remark}

\subsection{Other projective varieties}

\begin{theorem}%\label{theorem:general-base-obstruction}
Let $V$ be smooth, connected, and projective of dimension
$n\geq 3$.
For $d\geq 6$, set $X=\operatorname{SP}^{d}(V)$ and
$N=nd$.
Then there exists a nonalgebraic class
$\alpha_{V,d}\in H_{2N-4}(X,\mathbb{Z})$
of exact order two whose image generates a direct
summand $\mathbb{Z}/2$ of $D_{N-2}(X)$.
Furthermore, if $V$ admits an affine paving, each cycle class map of $X$
is split injective with finite cokernel annihilated by
$d!$, and
\[
H_{2N-4}(X,\mathbb{Z})
\cong
\operatorname{CH}_{N-2}(X)
\oplus\mathbb{Z}/2\oplus D'_{V,d},
\]
where $D'_{V,d}$ is finite and annihilated by $d!$.
In this case, the displayed Chow group and homology group
are not isomorphic as abelian groups.
\end{theorem}

\begin{proof}
The collision locus has codimension $n$, 
so localization yields the analogous commutative diagram to \eqref{equation:configuration-comparison}.
Three disjoint holomorphic coordinate balls and $d-6$ fixed points produce
the same map from $(\mathbb{RP}^{2n-1})^3$ and bundle restriction.
The Steenrod square calculation, the vanishing lemma, and the coefficient functional therefore apply without modification..

Under the additional affine-paving hypothesis,
by invoking
Theorem~\ref{theorem:affine-paving-comparison} and
Corollary~\ref{corollary:paving-bounds},
the torsion-cardinality argument establishes the final assertion.
\end{proof}

\begin{example}
For $X=\operatorname{SP}^{6}(\mathbb{P}^{3})$ , one has $N=18$ and
\[
H_{32}(X,\mathbb{Z})
\cong
\operatorname{CH}_{16}(X)
\oplus\mathbb{Z}/2\oplus D'_{3,6},
\qquad
720D'_{3,6}=0.
\]
The ranks of $\operatorname{CH}_i(X)$, for $0\leq i\leq18$,
are
\[
1,1,2,3,4,5,7,7,8,8,8,7,7,5,4,3,2,1,1.
\]
These are the coefficients of $u^i$ in
$[t^6]\prod_{j=0}^3(1-tu^j)^{-1}$ and sum to $84=\binom96$.
This information does not determine the full torsion subgroup or
the complementary group $D'_{3,6}$.
\end{example}

\section{Retractions and stabilization}

This section establishes integral stabilization under the addition of a point.
The stable cycle class map is split injective with finite cokernel annihilated by $i!$ in cycle dimension $i$.
It also provides an example showing that integral injectivity of multiplication cannot serve as a general step in a stability argument.

\subsection{Retaining a subconfiguration}
For this subsection let $V$ be a projective variety with a chosen point $o\in V$, 
and write$X_d(V)=\operatorname{SP}^{d}(V).$
For $0\leq r\leq s$, define
\[
\iota_{r,s}\colon X_r(V)\longrightarrow X_s(V),
\qquad c\longmapsto c+(s-r)o,
\]
and write $\iota_r=\iota_{r,r+1}$.
We adopt the conventions $\mathfrak S_0=\{1\}$ and $X_0(V)=\operatorname{Spec}(\mathbb{C})$.
For $0\leq r\leq d$, let $H_{d,r}=\mathfrak{S}_r\times\mathfrak{S}_{d-r}.$
The identification
$V^d/H_{d,r}=X_r(V)\times X_{d-r}(V)$ yields a proper first projection
$p_{d,r}$. Define the \emph{retention operators}
\begin{equation}\label{eq:retention}
 \begin{split}
 T^{\operatorname{CH}}_{d,r}&=(p_{d,r})_*\operatorname{tr}^{\mathfrak S_d}_{H_{d,r}}:
       \operatorname{CH}_i(X_d(V))\longrightarrow\operatorname{CH}_i(X_r(V)),\\
 T^{H}_{d,r}&=(p_{d,r})_*\operatorname{Tr}^{\mathfrak S_d}_{H_{d,r}}:
       H_k(X_d(V),\mathbb{Z})\longrightarrow H_k(X_r(V),\mathbb{Z}).
 \end{split}
\end{equation}
These maps preserve cycle dimension and homological
degree, respectively, and commute with cycle class maps.
We write $T_{d,r}$ when the same assertion applies to both theories.
Set $T_{d,d}=\operatorname{id}$, and interpret an operator with index outside $0\leq r\leq d$ as zero.

\begin{lemma}\label{lemma:retention-recurrence}
For $d\geq 1$ and $0\leq r\leq d-1$,   in both Chow groups and homology,
\[
T_{d,r}(\iota_{d-1})_*
=
T_{d-1,r}
+
(\iota_{r-1})_*T_{d-1,r-1}.
\]
The second term vanishes when $r=0$.
\end{lemma}

\begin{proof}
We first recall the double-coset description of the
transfer.
If $K\leq G$ and $Z$ is a projective $K$-variety, 
set $W=G\times_K Z$, the disjoint union of copies of $Z$ indexed by $G/K$. 
One has the identifications
\[
W/G=Z/K,
\qquad
W/H
=
\coprod_{[g]\in H\backslash G/K}
Z/(K\cap g^{-1}Hg).
\]
On the component indexed by $g$, the relative transfer from $W/G$
identifies with the transfer from $Z/K$ to
$Z/(K\cap g^{-1}Hg)$, followed by the indicated conjugation identification.
For algebraic cycles this follows directly from \eqref{equation:relative-transfer}:
the stabilizer of $[g,z]$ in $G$ is $gK_zg^{-1}$, and its stabilizer in $H$
is $H\cap gK_zg^{-1}$. The same quotient of stabilizer orders therefore
occurs in the two formulas. For homology, decompose the coset sum
\eqref{equation:coset-transfer} according to the same double cosets. This proves the
description integrally in both theories.

Now let $G=\mathfrak{S}_d, K=\mathfrak{S}_{d-1},$ and $Z=V^{d-1}$, 
where $K$ fixes the final coordinate label.
There is a proper $G$-equivariant map
\[
 f:G\times_K V^{d-1}\longrightarrow V^d,
 \qquad [g,z]\longmapsto g(z,o).
\]
The quotient map $f_G$ coincides with $\iota_{d-1}$.
For $H=H_{d,r}$, the double cosets are classified by whether the distinguished
last label lies in the retained block of size $r$ or to its complement.
If it is discarded, the component transfer, followed by $(p_{d,r}f_H)_*$, is $T_{d-1,r}$. 
If it is retained, it recovers $(\iota_{r-1})_*T_{d-1,r-1}$. 
For $r=0$ only the first case arises.
Proposition~\ref{proposition:transfer-naturality}
therefore yields the formula in both theories.

In particular, the argument does not require the points of $z$ to be distinct or distinct from $o$. 
To verify the multiplicities, 
consider a zero-cycle \(c=\sum_x m_x[x]\): 
the coefficient associated to the retained subcycle \(\sum_x b_x[x]\) is $\prod_x\binom{m_x}{b_x}$. 
Meanwhile, at $o$, the two cases above correspond exactly to the identity
$\binom{m_o+1}{b_o}=\binom{m_o}{b_o}+\binom{m_o}{b_o-1}.$
The double coset argument, not merely this pointwise verification, establishes the identity for all cycles and homology classes.
\end{proof}

\begin{proposition}\label{proposition:integral-retraction}
Given any projective $V$ with a chosen point, all $d\geq1$ and $i,k\geq0$, 
the maps induced by $\iota_{d-1}$ on Chow groups and homology admit integral left inverses
\[
R^{\bullet}_d
=
\sum_{r=0}^{d-1}
(-1)^{d-1-r}
(\iota_{r,d-1})_*T^{\bullet}_{d,r},
\qquad
\bullet\in\{\operatorname{CH},H\}.
\]
They satisfy
\[
R^{\bullet}_d(\iota_{d-1})_*=\operatorname{id},
\qquad
R^H_d\operatorname{cl}_{X_d(V),i}
=
\operatorname{cl}_{X_{d-1}(V),i}R^{\operatorname{CH}}_d.
\]
Consequently, the map $D_i(X_{d-1}(V))\longrightarrow D_i(X_d(V))$ is split injective.
\end{proposition}

\begin{proof}
Substituting  the recurrence from
Lemma~\ref{lemma:retention-recurrence}, one finds that the term indexed by $r$ in the first sum cancels the
term indexed by $r+1$ in the second sum.
This means that only $T_{d-1,d-1}=\operatorname{id}$ survives.
Compatibility with cycle class maps holds termwise from \eqref{eq:retention}. 
It also yields the induced retraction $R_d^H$ on cokernels,
for which the same left inverse identity holds.
\end{proof}

For curves, fixed-point embeddings admit a classical Chow-theoretic treatment in \cite{Collino}.
A higher-dimensional injectivity result with rational coefficients appears in \cite[Thm. 2.3]{Banerjee}.
The proposition above constructs integral retractions and,
critically for later applications, compatible retractions on homology and on the cokernels of cycle class maps.

\subsection{The stable range}

Return to $V=\mathbb{P}^{n}$.
Set
\[
Y_n=\operatorname{SP}^{\infty}(\mathbb{P}^{n})
=
\varinjlim_d X_{n,d},
\qquad
\operatorname{CH}_q(Y_n)
=
\varinjlim_d\operatorname{CH}_q(X_{n,d}),
\]
where the transition maps are given by adding the point $o$.
We equip $Y_n$ with the weak topology.
Singular homology commutes with this direct limit; 
equivalently, this follows from the compatible CW filtrations of finite symmetric products \cite[\S4.K, P.~483]{HatcherTopology}.
Thus the compatible finite-stage cycle maps define
\[
\operatorname{cl}_{\infty,q}\colon
\operatorname{CH}_q(Y_n)
\longrightarrow H_{2q}(Y_n,\mathbb{Z}).
\]
Write $D_{\infty,q}(n)=\operatorname{coker}(\operatorname{cl}_{\infty,q}).$

\begin{theorem}\label{theorem:stable-range}
For all $n,q\geq0$ and $d\geq q$, addition of $o$ induces an isomorphism
\[
(\iota_d)_*\colon
\operatorname{CH}_q(X_{n,d})
\longrightarrow
\operatorname{CH}_q(X_{n,d+1}).
\]
The stable cycle class map is split injective,  
and there exist natural identifications
\[
\operatorname{CH}_q(Y_n)
\cong\operatorname{CH}_q(X_{n,q}),
\qquad
D_{\infty,q}(n)\cong D_q(X_{n,q}).
\]
In particular, $D_{\infty,q}(n)$ is finite and annihilated by $q!$.
For $d\geq q$, the torsion of $\operatorname{CH}_q(X_{n,d})$ and $H_{2q}(X_{n,d},\mathbb{Z})$,
together with $D_q(X_{n,d})$, is annihilated by $q!$.
The stable rank is
\begin{equation}\label{eq:stable-rank}
\operatorname{rank}\operatorname{CH}_q(Y_n)=[z^q]\prod_{j=1}^{n}(1-z^j)^{-1}.
\end{equation}
For $n=0$, the product is interpreted as $1$.
\end{theorem}

\begin{proof}
The case $n=0$ is immediate.
For $n\geq1$, since projective space is simply connected,
the connectivity result for symmetric products
\cite[Thm. 1.3]{Kallel} yields
\[
(\iota_d)_*\colon
H_{2q}(X_{n,d},\mathbb{Z})
\xrightarrow{\ \cong\ }
H_{2q}(X_{n,d+1},\mathbb{Z})
\qquad(d\geq q).
\]
Hence its left inverse $R^H_{d+1}$ from Proposition~\ref{proposition:integral-retraction} is also a right inverse.
For $\zeta\in\operatorname{CH}_q(X_{n,d+1})$,
compatibility yields
\[
\operatorname{cl}_{X_{n,d+1},q}
\bigl(\zeta-(\iota_d)_*R^{\operatorname{CH}}_{d+1}(\zeta)\bigr)
=0.
\]
The cycle class map $\operatorname{cl}_{X_{n,d+1},q}$ is injective by
Theorem~\ref{theorem:main-comparison},
so $\zeta-(\iota_d)_*R^{\operatorname{CH}}_{d+1}(\zeta)=0$.
Therefore $(\iota_d)_*$ is surjective on Chow groups.
It is injective by Proposition~\ref{proposition:integral-retraction}.
This proves $(\iota_d)_*$ is an isomorphism.

Observe that both direct systems are constant up to canonical isomorphism starting at stage $q$. 
Their cycle class maps, cokernels, and torsion groups are therefore identified with those at the corresponding finite stage.
Finally, for $d\geq q$, the tuples in \eqref{equation:rank-formula} correspond to partitions
of $q$ with parts at most $n$, padded with zeros. This establishes the rank formula \eqref{eq:stable-rank}.
\end{proof}

\subsection{Persistence of the obstruction}

\begin{theorem}\label{theorem:persistence}
Let $n\geq3$, $m\geq6$, and $q=nm-2$.
For each $d\geq m$, define
\[
\alpha_{n;m,d}
=
(\iota_{m,d})_*\alpha_{n,m}
\in H_{2q}(X_{n,d},\mathbb{Z}).
\]
This class has exact order two and is nonalgebraic,
and its image generates a direct summand $\mathbb{Z}/2$ of
$D_q(X_{n,d})$ and of $D_{\infty,q}(n)$.
In particular,
\[
H_{2q}(X_{n,d},\mathbb{Z})
\cong
\operatorname{CH}_q(X_{n,d})
\oplus\mathbb{Z}/2\oplus D''_{n;m,d},
\]
where $D''_{n;m,d}$ is finite.
Thus the Chow group and homology group are not isomorphic as abstract abelian groups.
\end{theorem}

\begin{proof}
Iterating
Proposition~\ref{proposition:integral-retraction} produces split injections on homology and cycle class
cokernels from stage $m$ to each later stage.
Thus the order two class remains nonzero.
The retraction on cokernels, composed with the detector $\overline{\chi}$ at stage $m$, 
yields a homomorphism to $\mathbb{Z}/2$ sending its image to $1$.
This establishes nonalgebraicity together with the direct summand claim.

Since $q=nm-2\geq m$, the stable cokernel is identified
with the cokernel at stage $q$ by
Theorem~\ref{theorem:stable-range}.
The remaining assertions follow from
Theorem~\ref{theorem:main-comparison} and the
torsion-cardinality argument.
\end{proof}

\begin{proof}[Proof of Theorem~\ref{theorem:main-stability}]
This follows by combining Theorems~\ref{theorem:stable-range} and~\ref{theorem:persistence}, 
taking $(n,m,q)=(3,6,16)$ for the explicit example displayed above.
\end{proof}

\begin{example}
Take $(n,m,q)=(3,6,16)$.
At stage $d=16$, both the Chow group and the homology
group in degree $32$ have rank
\[
[z^{16}]
\frac{1}{(1-z)(1-z^2)(1-z^3)}
=
\sum_{c=0}^{5}
\left(
\left\lfloor\frac{16-3c}{2}\right\rfloor+1
\right)
=
30.
\]
They are nonisomorphic, since the homology cokernel
contains the propagated direct summand $\mathbb{Z}/2$.
\end{example}

\subsection{Multiplication of effective cycles}

For $e\geq1$, define
\[
\mu_{e,d}\colon X_{n,d}\longrightarrow X_{n,ed},
\qquad c\longmapsto ec.
\]
Its action on Chow groups is given by proper pushforward along this morphism.

\begin{example}%\label{example:projective-line-multiplication}
Under the identification
$\operatorname{SP}^{d}(\mathbb{P}^{1})\cong\mathbb{P}^{d}$,
the map $\mu_{e,d}$ sends a binary form $F$ to $F^e$.
Thus $\mu_{e,d}^{*}\mathcal{O}_{\mathbb{P}^{ed}}(1)=\mathcal{O}_{\mathbb{P}^{d}}(e).$
The projection formula yields
\[
(\mu_{e,d})_*[\mathbb{P}^{q}]
=
e^q[\mathbb{P}^{q}]
\quad\text{in}\quad
\operatorname{CH}_q(\mathbb{P}^{ed}),
\qquad 0\leq q\leq d.
\]
In particular, for $e>1$ and $q>0$, this is not an
integral isomorphism.
For $q=2$, it is not scalar multiplication by the scalar $e$.
By contrast, $\iota_d:\mathbb{P}^d\hookrightarrow\mathbb{P}^{d+1}$
is a linear embedding and carries each standard Chow generator to the corresponding standard generator.
\end{example}

\begin{proposition}%\label{proposition:multiplication-torsion}
For all $d\geq2, e\geq1$, there exists class $\tau_d\in\operatorname{CH}_3(X_{2,d})$ of exact order two such that
$(\mu_{e,d})_*\tau_d=e\tau_{ed}$.
Consequently, every even multiplicity annihilates $\tau_d$, 
whereas every odd multiplicity preserves its nonzero stabilized class.
In particular, $\mu_{2,3}:X_{2,3}\to X_{2,6}$ is not injective on $\operatorname{CH}_3$, 
and multiplication need not induce an integral Chow injection even in the range $q\leq d$.
\end{proposition}

\begin{proof}
Let $B\subset X_{2,2}$ be the image of the diagonal, and
set
\[
F=(\mathbb{P}^{2})^2\setminus\Delta,
\qquad
U=X_{2,2}\setminus B=F/\mathfrak{S}_2.
\]
The double cover $F\longrightarrow U$ together with the sign character 
define an algebraic line bundle $\mathcal{L}$ satisfying $\mathcal{L}^{\otimes2}\cong\mathcal{O}_U$.
The complement of the diagonal in $(\mathbb{P}^{2})^2$ is simply connected.
Indeed, general position allows loops and their bounding
disks to avoid a smooth submanifold of real codimension four.
Hence
\[
\pi_1(U)\cong\mathbb{Z}/2,
\qquad
H_1(U,\mathbb{Z})\cong\mathbb{Z}/2,
\qquad
H^1(U,\mathbb{Z})=0.
\]
If $w\in H^1(U,\mathbb{F}_2)$ classifies this cover, 
then $c_1(\mathcal{L})=\beta_2(w).$
Since $w\neq0$ and $H^1(U, \mathbb{Z})=0$, the coefficient sequence forces this Bockstein class nonzero.
Thus $\mathcal{L}$ has exact order two in $\operatorname{Pic}(U)$, and its topological first Chern class is nonzero.

Since $\dim X_{2,2}=4$ and $\dim B=2$, localization yields
compatible isomorphisms
\[
\operatorname{CH}_3(X_{2,2})
\cong
\operatorname{CH}^{1}(U)
=
\operatorname{Pic}(U),
\]
\[
H_6(X_{2,2},\mathbb{Z})
\cong H^2(U,\mathbb{Z}).
\]
Let $\tau_2$ be the class corresponding to $\mathcal{L}$, and define
$\tau_d=(\iota_{2,d})_*\tau_2$.
The Steenrod--Nakaoka injection preserves its nonzero cycle class,
so every $\tau_d$ has exact order two.

Let $Y=Y_2$.
The Dold--Thom theorem and its product consequence yield
\[
Y\simeq K(\mathbb{Z},2)\times K(\mathbb{Z},4);
\]
see \cite[Thm. 4K.6 and Cor. 4K.7]{HatcherTopology}.
Here $K(A,r)$ denotes an Eilenberg--Mac Lane space with $\pi_r=A$ and all other positive homotopy groups trivial.
Thus, in low degrees,
$H_1(Y,\mathbb{Z})=H_3(Y,\mathbb{Z})=H_5(Y,\mathbb{Z})=0,$
while $H_2(Y,\mathbb{Z})$ and $H_4(Y,\mathbb{Z})$ are free abelian.
For this degree-five vanishing assertion, one may construct $K(\mathbb{Z},4)$ from $S^4$ by attaching cells of
dimension at least six; the resulting cellular chain group in degree five vanishes.

Let $t\in H_6(Y,\mathbb{Z})$ denote the stabilized cycle class of $\tau_d$;
it is independent of $d$ and has exact order two.
In degree six, the integral K\"unneth formula has no
$\operatorname{Tor}$ summands.
Aside from the two edge summands, its only nonzero summands are
$H_2(Y,\mathbb{Z})\otimes H_4(Y,\mathbb{Z})$ and
$H_4(Y,\mathbb{Z})\otimes H_2(Y,\mathbb{Z})$, both of which are free.
Consequently, the diagonal $\delta:Y\to Y\times Y$ satisfies
$\delta_*(t)=t\otimes1+1\otimes t.$
Indeed, mixed components of a torsion element vanish, and projection
onto either factor identifies each edge component with $t$.
Addition on $Y$ therefore sends the diagonal image to $2t$,
and the iterated diagonal carries $t$ to the sum of its $e$ edge terms.
Thus multiplication by $e$ on the symmetric product monoid sends $t$
to $et$, which vanishes for even $e$.
Therefore, 
addition in the infinite symmetric product 
sends the iterated diagonal of $t$ to $et$
since the finite maps satisfy
\[
\mu_{e,d+1}\circ \iota_d
=
\iota_{ed,e(d+1)}\circ\mu_{e,d},
\]
and they therefore induce this multiplication map on the direct limit.
Since $H_6(X_{2,ed},\mathbb{Z})\to H_6(Y,\mathbb{Z})$ is injective, the equality
between the stabilized cycle classes of $(\mu_{e,d})_*\tau_d$ and
$e\tau_{ed}$ already holds at the finite stage $ed$.
The target cycle map is injective by Theorem~\ref{theorem:main-comparison}, proving
$(\mu_{e,d})_*\tau_d=e\tau_{ed}$.
\end{proof}

\begin{remark}
The preceding proposition concerns multiplication, not
addition of a fixed point.
It therefore does not contradict integral stabilization result in Theorem~\ref{theorem:stable-range}.

For the Chow variety
$C_{p,d}(\mathbb{P}^{n})$ of effective $p$-dimensional
cycles of degree $d$ in $\mathbb{P}^{n}$, one may similarly define
\[
\mathcal{C}_p(\mathbb{P}^{n})
=
\varinjlim_d C_{p,d}(\mathbb{P}^{n})
\]
by addition of a fixed $p$-plane, and define its Chow groups via this direct system.
Lawson's complex suspension theorem \cite[Thm. 3]{Lawson} addresses the topology of these completed spaces.
A comparison for their Chow groups requires a separate algebraic argument.
The present paper do not establish integral Chow suspension or finite-degree Chow stability for $p>0$.
In the zero-cycle parameter spaces considered here,
Theorem~\ref{theorem:persistence} instead proves stable nonsurjectivity
in the specified degrees.
\end{remark}

\section{Unramified cohomology}

This section is devoted to connecting the codimension-two cokernel to unramified cohomology
on the smooth configuration open set,
following Colliot-Thélène–Voisin \cite[Thm.~3.7]{ColliotTheleneVoisin}.
It also clarifies that this obstruction cannot produce an integral Hodge counterexample
on a smooth projective resolution,
since these resolutions are rational and unramified classes on the open set must ramify along the boundary.

\subsection{A comparison for smooth open varieties}

Assume that $M$ is a smooth connected complex quasiprojective variety.
Set
\[
D^2_{\mathrm{coh}}(M)
=
H^4(M,\mathbb{Z}(2))/
\operatorname{im}(\operatorname{cl}^{2}_{M}),
\]
where
$\mathbb{Z}(r)=(2\pi\sqrt{-1})^r\mathbb{Z}$.
For $s\geq1$, let
$\mathbb{Z}/s(r)=\mathbb{Z}(r)/s\mathbb{Z}(r)$.
After fixing these generators, the underlying groups
are identified with the untwisted groups used earlier.
For a coefficient group $A$, let
$\mathcal{H}^{b}_{M}(A)$ denote the Zariski sheaf associated
to $W\longmapsto H^b(W,A)$, and define
\[
H^b_{\mathrm{nr}}(M,A)
=
H^0(M,\mathcal{H}^{b}_{M}(A)).
\]
For an abelian group $A$, write
$A[s]=\ker(s\colon A\longrightarrow A)$.

We require the following form of the Colliot-Th\'el\`ene--Voisin argument.
Its target is the full cohomological cokernel; 
we make no assumption of a pure Hodge structure on a nonproper variety.

\begin{proposition}\label{proposition:unramified-sequence}
For each smooth connected complex quasiprojective $M$ and integer $s\geq1$, 
there exists a natural exact sequence
\[
0\longrightarrow
H^3_{\mathrm{nr}}(M,\mathbb{Z}(2))/s
\longrightarrow
H^3_{\mathrm{nr}}(M,\mathbb{Z}/s(2))
\longrightarrow
D^2_{\mathrm{coh}}(M)[s]
\longrightarrow0.
\]
Passing to the direct limit yields
\[
0\longrightarrow
H^3_{\mathrm{nr}}(M,\mathbb{Z}(2))
\otimes\mathbb{Q}/\mathbb{Z}
\longrightarrow
H^3_{\mathrm{nr}}(M,\mathbb{Q}/\mathbb{Z}(2))
\longrightarrow
D^2_{\mathrm{coh}}(M)_{\mathrm{tors}}
\longrightarrow0.
\]
If $H^3(M,\mathbb{Q})=0$ and
$\operatorname{cl}^{2}_{M}$ is injective, then
$H^3_{\mathrm{nr}}(M,\mathbb{Z}(2))=0.$
\end{proposition}

\begin{proof}
The Gersten resolution and coniveau spectral sequence for smooth varieties give
\[
E_2^{a,b}
=
H^a(M,\mathcal{H}^{b}_{M}(\mathbb{Z}(2)))
\Longrightarrow H^{a+b}(M,\mathbb{Z}(2)),
\qquad
E_2^{a,b}=0\quad(a>b).
\]
See \cite[Thm.~ 4.2]{BlochOgus} and
\cite[Thm.~ 2.3 and Cor.~ 2.5]{ColliotTheleneVoisin}.
Set $C^2(M)=H^2(M,\mathcal{H}^{2}_{M}(\mathbb{Z}(2))).$
The codimension-two term of the Gersten complex is the group of codimension-two cycles; 
the boundary operator attached to a rational function is its divisor.
Consequently the cycle map factors as
\[
\operatorname{CH}^{2}(M)
\twoheadrightarrow C^2(M)
\longrightarrow H^4(M,\mathbb{Z}(2))
\]
where image of the last arrow is the algebraic subgroup.
Only the surjection onto $C^2(M)$ is used; we do not identify this group
with $\operatorname{CH}^2(M)$ in general.
Therefore, the filtration in total degree four gives an exact sequence
\begin{equation}\label{equation:coniveau-cokernel}
0\longrightarrow
H^1(M,\mathcal{H}^{3}_{M}(\mathbb{Z}(2)))
\longrightarrow D^2_{\mathrm{coh}}(M)
\longrightarrow
H^0(M,\mathcal{H}^{4}_{M}(\mathbb{Z}(2))).
\end{equation}
Indeed, $E_{\infty}^{1,3}=E_2^{1,3}$, 
and the algebraic subgroup is the second step of the coniveau filtration in
degree four.

By \cite[Thm.~ 3.1]{ColliotTheleneVoisin}, 
the integral cohomology sheaves $\mathcal H_M^b(\mathbb{Z}(2))$ are torsion free, 
and there exist short exact sequences
\[
0\longrightarrow
\mathcal{H}^{b}_{M}(\mathbb{Z}(2))
\xrightarrow{\ s\ }
\mathcal{H}^{b}_{M}(\mathbb{Z}(2))
\longrightarrow
\mathcal{H}^{b}_{M}(\mathbb{Z}/s(2))
\longrightarrow0.
\]
The final group in \eqref{equation:coniveau-cokernel} is torsion free.
Thus
\[
H^1(M,\mathcal{H}^{3}_{M}(\mathbb{Z}(2)))[s]
\cong D^2_{\mathrm{coh}}(M)[s].
\]
The long exact sequence of Zariski sheaf cohomology for $b=3$ proves the first assertion.
Compatible coefficient inclusions and filtered direct limits yield the second one, 
as in \cite[proof of Thm~ 3.7]{ColliotTheleneVoisin}.
The argument underlying that theorem therefore applies to our full cokernel without imposing
Hodge classes on $M$.

For the last assertion, the spectral sequence also gives
\[
H^3(M,\mathbb{Z}(2))
\longrightarrow
H^3_{\mathrm{nr}}(M,\mathbb{Z}(2))
\longrightarrow C^2(M)
\longrightarrow H^4(M,\mathbb{Z}(2)).
\]
If the composite map from $\operatorname{CH}^2(M)$ to $H^4(M,\mathbb{Z}(2))$ is injective,
its surjective factorization through $C^2(M)$ shows that
$C^2(M)\to H^4(M,\mathbb{Z}(2))$ is injective.
Thus the map $H^3(M,\mathbb{Z}(2)) \to H^3_{\mathrm{nr}}(M,\mathbb{Z}(2))$ is surjective.
Its source is torsion whenever $H^3(M,\mathbb{Q})=0$, whereas its target is
torsion free by \cite[Thm.~3.1]{ColliotTheleneVoisin}. 
Therefore $H^3_{\mathrm{nr}}(M,\mathbb{Z}(2))$ vanishes.
\end{proof}

\begin{theorem}\label{theorem:configuration-unramified}
Let $n\geq3$ and $d\geq1$.
Then $H^3_{\mathrm{nr}}(U_{n,d},\mathbb{Z}(2))=0$,
and there exist natural isomorphisms
\[
H^3_{\mathrm{nr}}(U_{n,d},\mathbb{Z}/s(2))
\cong D_{nd-2}(X_{n,d})[s],
\]
\[
H^3_{\mathrm{nr}}(U_{n,d},\mathbb{Q}/\mathbb{Z}(2))
\cong D_{nd-2}(X_{n,d}).
\]
The latter group is finite and annihilated by $d!$.
For $d\geq6$, it contains a direct summand
$\mathbb{Z}/2$.
\end{theorem}

\begin{proof}
Removing the collision diagonals of codimension $n\geq3$ from $(\mathbb{P}^n)^d$ preserves rational cohomology in degree three.
Consequently, $H^3(F_d(\mathbb{P}^{n}),\mathbb{Q})=0.$
The finite covering $F_d(\mathbb{P}^n)\to U_{n,d}$ admits a transfer map;
this makes pullback from $U_{n,d}$ in rational cohomology injective, 
so $H^3(U_{n,d},\mathbb{Q})=0$ as well.
The localization comparison
\eqref{equation:configuration-comparison} holds for all $d\geq1$.
Theorem~\ref{theorem:main-comparison} therefore yields
injectivity of $\operatorname{cl}^{2}_{U_{n,d}}$ and
\[
D^2_{\mathrm{coh}}(U_{n,d})
\cong D_{nd-2}(X_{n,d}).
\]
The theorem now follows by applying Proposition~\ref{proposition:unramified-sequence},
Theorems~\ref{theorem:main-comparison} and Theorems~\ref{theorem:main-obstruction}.
\end{proof}

\subsection{Smooth projective resolutions}

For a smooth projective variety $Y$, let
\[
Z^4(Y)
=
\operatorname{Hdg}^{4}(Y,\mathbb{Z}(2))/
\operatorname{im}(\operatorname{cl}^{2}_{Y}).
\]
All torsion classes are included among the integral
Hodge classes.
Since $H^4(Y,\mathbb{Z}(2))/\operatorname{Hdg}^{4}(Y,\mathbb{Z}(2))$ is torsion free, 
we have
\[
D^2_{\mathrm{coh}}(Y)_{\mathrm{tors}}
=
Z^4(Y)_{\mathrm{tors}}.
\]
Theorem~3.7 of \cite{ColliotTheleneVoisin} identifies this torsion group as the quotient
of $H^3_{\mathrm{nr}}(Y,\mathbb{Q}/\mathbb{Z}(2))$ by the image of
$H^3_{\mathrm{nr}}(Y,\mathbb{Z}(2))\otimes\mathbb{Q}/\mathbb{Z}$.
Theorem~\ref{theorem:configuration-unramified} adapts this same method
to a nonproper smooth configuration space and its full cycle-class cokernel; 
it is not a new proof of the cited projective theorem.

Voisin's vanishing theorem for uniruled threefolds and for the
Calabi--Yau threefolds specified in \cite[Thm.~2]{Voisin}
concerns $Z^4$ on smooth projective threefolds.
Our $\operatorname{SP}^6(\mathbb{P}^3)$ is singular and of dimension $18$, not $3$.
More decisively, its smooth projective models have vanishing integral Hodge defect, as shown below.

\begin{proposition}%\label{proposition:rational-resolution}
Let $n\geq3$, $d\geq1$, $X=X_{n,d}$, and $U=U_{n,d}$.
The variety $X$ is rational.
Let $\pi\colon\widetilde{X}\longrightarrow X$ be a projective log resolution 
which is an isomorphism over $U$, and let $E_1,\ldots,E_t$ denote the irreducible divisors of $\widetilde{X}\setminus U$.
Then
\[
Z^4(\widetilde{X})=0,
\qquad
H^3_{\mathrm{nr}}
(\widetilde{X},\mathbb{Z}/s(2))=0
\quad(s\geq1).
\]
The boundary residues induce an injection
\[
D_{nd-2}(X)[s]
\longrightarrow
\bigoplus_{j=1}^{t}
H^2(\mathbb{C}(E_j),\mathbb{Z}/s(1)).
\]
Here cohomology of a function field means Galois cohomology.
For $d\geq6$, the class
$a=\beta_2(w_3(E))\in H^4(U,\mathbb{Z})$
does not extend to an integral Hodge class on
$\widetilde{X}$.
\end{proposition}

\begin{proof}
Rationality is classical; see \cite{Mattuck}.
For completeness, consider a general unordered configuration in an affine chart
given by $(t_i,u_{i,2},\ldots,u_{i,n}), 1\leq i\leq d$, with pairwise distinct $t_i$.
It determines
$P(T)=\prod_{i=1}^{d}(T-t_i)$
and interpolation polynomials $Q_j$ satisfying
$\deg Q_j<d, Q_j(t_i)=u_{i,j}$ for $2\leq j\leq n.$
The $d$ nonleading coefficients of $P$ and the
$d(n-1)$ coefficients of the $Q_j$ provide birational coordinates.
Conversely, a separable polynomial $P$ together with the $Q_j$ recover the unordered configuration.
Thus $X$ and $\widetilde{X}$ are rational.

The group $Z^4$ is a birational invariant of smooth projective varieties; see \cite[Introduction]{Voisin}.
The blowup formula adds only degree-two Hodge classes and degree-zero classes from the center, 
which are algebraic by the integral Lefschetz $(1,1)$ theorem.
Weak factorization therefore implies the birational invariance. 
Hence $Z^4(\widetilde{X})=Z^4(\mathbb{P}^{nd})=0.$
Unramified cohomology in positive degree with these finite coefficients vanishes for smooth projective
rational varieties; see
\cite[Thm. 2.8]{ColliotTheleneVoisin}.
This proves the second vanishing statement.

By the Gersten description \cite[Thm.~2.3]{ColliotTheleneVoisin}, 
a class in $H^3_{\mathrm{nr}}(U,\mathbb{Z}/s(2))$
is a function-field class unramified at all codimension
one points of $U$.
If its residues along each $E_j$ also vanish, 
it is unramified at all codimension one points of $\widetilde X$.
This means it lies in the trivial group $H^3_{\mathrm{nr}}(\widetilde X,\mathbb{Z}/s(2))$.
Theorem~\ref{theorem:configuration-unramified} now yields
the injection.
In particular, the nonzero classes detected in
Theorem~\ref{theorem:configuration-unramified} must ramify along some boundary divisor.

Finally, suppose $a$ admits an extension to an integral Hodge class on $\widetilde X$.
Since $Z^4(\widetilde{X})=0$, this extended class would be algebraic.
Restricting the representing Chow class to $U$ would force $a$ to be algebraic, 
contradicting Proposition~\ref{proposition:explicit-obstruction}.
\end{proof}

Thus the finite cokernel for the singular symmetric
product is not the integral Hodge defect of a smooth
projective resolution.
The vanishing of $Z^4(\widetilde{X})$ does not assert that
each degree-four cohomology class is algebraic, since
non-Hodge classes are not counted in $Z^4$.

\subsection{Final remark}
For $d=0$ or $n=0$, $X_{n,d}$ is a point; 
for $d=1$ its cycle class map is an isomorphism. 
The same holds for $n=1$ since $\operatorname{SP}^d(\mathbb{P}^1)\cong \mathbb{P}^d$.
The explicit construction in Section~\ref{sec:obstruction} requires
three moving pairs and $n\geq3$.
We do not claim that $(n,d)=(3,6)$ is minimal.
This result does not determine the complementary finite cokernels,
their optimal annihilators, or all remaining cases with $n\leq2$ or
$d\leq5$.
The persistence theorem does, however, disprove integral surjectivity
both in some finite stages with $q\leq d$ and in the infinite symmetric
product. Integral Chow suspension for general positive-dimensional Chow varieties is not asserted.

\appendix

\section{Codimension two injectivity from a finite \'etale cover}

The argument in this appendix is independent of the coarse quotient transfers and Sylow descent.
It proves codimension two injectivity and the factorial bound on the cokernel, but does not provide the general splitting.

\begin{lemma}\label{lemma:torsion-injectivity}
Let $M$ be smooth and quasiprojective with
$H^3(M,\mathbb{Q})=0$.
Then
\[
\operatorname{cl}^{2}_{M}\colon
\operatorname{CH}^{2}(M)_{\mathrm{tors}}
\longrightarrow H^4(M,\mathbb{Z})_{\mathrm{tors}}
\]
is injective.
\end{lemma}

\begin{proof}
We apply the diagram for smooth varieties in
\cite[proof of Cor.~3.5]{Totaro}, which follows from
Bloch--Merkurjev--Suslin theory.
For a prime $\ell$, set
\[
B_{\ell}(M)
=
H^1_{\mathrm{Zar}}
\bigl(
M,\mathcal{H}^{2}_{\mathrm{\acute{e}t}}
(\mathbb{Q}_{\ell}/\mathbb{Z}_{\ell}(2))
\bigr)
\]
where $\mathcal H_{\mathrm{et}}^{\,2}$ denotes the Zariski sheaf associated to second \'etale cohomology, 
and $(2)$ denotes the second Tate twist.
Then there exists a commutative diagram
\[
\begin{tikzcd}
B_{\ell}(M)
  \arrow[r,hook,"\gamma_{\ell}"]
  \arrow[d,two heads,"b_{\ell}"']
&
H^3_{\mathrm{\acute{e}t}}
(M,\mathbb{Q}_{\ell}/\mathbb{Z}_{\ell}(2))
  \arrow[d,"\delta_{\ell}"]
\\
\operatorname{CH}^{2}(M)[\ell^{\infty}]
  \arrow[r,"\operatorname{cl}_{\ell}"']
&
H^4_{\mathrm{\acute{e}t}}
(M,\mathbb{Z}_{\ell}(2)),
\end{tikzcd}
\]
where $b_{\ell}$ is surjective and $\gamma_{\ell}$ is
injective.
Here $[\ell^{\infty}]$ denotes the $\ell$-primary
torsion subgroup, and $\delta_{\ell}$ is the coefficient
connecting homomorphism associated to
$0\to\mathbb{Z}_\ell(2)\to\mathbb{Q}_\ell(2)\to\mathbb{Q}_\ell/\mathbb{Z}_\ell(2)\to 0$.

Comparison with singular cohomology together with $H^3(M,\mathbb{Q})=0$ yield
$H^3_{\mathrm{\acute{e}t}}(M,\mathbb{Q}_{\ell}(2))=0.$
Hence $\delta_{\ell}$ is injective.
If a class $u\in\operatorname{CH}^{2}(M)[\ell^{\infty}]$ has
zero integral cycle class, choose $v\in B_{\ell}(M)$
with $b_{\ell}(v)=u$.
Then
$\delta_{\ell}\gamma_{\ell}(v)=0$.
Injectivity of both maps yields $v=0$ and $u=0$.
Applying this to each primary component of a torsion class gives the injectivity of $\operatorname{cl}^{2}_{M}$.
\end{proof}

\begin{proposition}%\label{proposition:independent-codimension-two}
Let $n\geq3$, $d\geq1$, $N=nd$, and $X=X_{n,d}$.
Then
\[
\operatorname{cl}_{X,N-2}\colon
\operatorname{CH}_{N-2}(X)
\longrightarrow H_{2N-4}(X,\mathbb{Z})
\]
is injective, with finite cokernel annihilated by $d!$,
independently of the arguments in Sections \ref{sec:transfers} and \ref{sec:splitting}.
\end{proposition}

\begin{proof}
Let $Y=(\mathbb{P}^{n})^d,
F=F_d(\mathbb{P}^{n}),
U=F/\mathfrak{S}_d,$
and let $f\colon F\longrightarrow U$.
Localization yields
\[
\operatorname{CH}^{2}(Y)
\cong\operatorname{CH}^{2}(F),
\qquad
H^4(Y,\mathbb{Z})\cong H^4(F,\mathbb{Z}).
\]
The product affine paving therefore makes
$\operatorname{cl}^{2}_{F}$ an isomorphism.
The same localization argument yields
$H^3(F,\mathbb{Q})=0$.
Finite cover transfer implies $f^\ast:H^3(U,\mathbb{Q})\to H^3(F,\mathbb{Q})$ injective,
so $H^3(U,\mathbb{Q})=0$.

If $\zeta\in\operatorname{CH}^2(U)$ has zero cycle class, 
injectivity on the covering space gives $f^\ast(\zeta)=0$. 
Flat pullback and proper pushforward yield $d!\zeta=f_\ast f^\ast(\zeta)=0$.
The kernel is therefore torsion, and Lemma~\ref{lemma:torsion-injectivity} shows it is zero.

For $z\in H^4(U,\mathbb{Z})$, choose
$\xi\in\operatorname{CH}^{2}(F)$ with
$\operatorname{cl}^{2}_{F}(\xi)=f^*(z)$.
The cohomological proper pushforward for this finite \'etale map, 
defined via the complex orientations, agrees with covering transfer.
Consequently,
\[
 \operatorname{cl}_U^2(f_\ast\xi)=f_\ast f^\ast (z)=d!z.
\]
Thus the cokernel is annihilated by $d!$.
The smooth quasiprojective varieties have finite CW homotopy type, by
semialgebraic triangulation of a compactification with its boundary.
Thus the cokernel is finite.
Finally, the localization and duality diagram \eqref{equation:configuration-comparison}, 
valid here even for $1\leq d<6$,
carries these assertions over from $U$ to $X$.
\end{proof}

For $d\geq6$, this proposition and
Proposition~\ref{proposition:explicit-obstruction}
also establish abstract nonisomorphism in codimension two
without invoking the splitting theorem.
The source embeds into a finitely generated homology group.
Its torsion image is a proper subgroup of the torsion
subgroup of the target because it does not contain
$\alpha_{n,d}$.

\end{document}